\documentclass[12pt, reqno]{amsart}
\usepackage{amsmath, amsthm, amscd, amsfonts, amssymb, graphicx, color}
\usepackage[bookmarksnumbered, colorlinks, plainpages]{hyperref}
\hypersetup{colorlinks=true,linkcolor=red, anchorcolor=green, citecolor=cyan, urlcolor=red, filecolor=magenta, pdftoolbar=true}

\newtheorem{theorem}{Theorem}[section]
\newtheorem{lemma}[theorem]{Lemma}
\newtheorem{proposition}[theorem]{Proposition}
\newtheorem{corollary}[theorem]{Corollary}
\theoremstyle{definition}
\newtheorem{definition}[theorem]{Definition}
\newtheorem{example}[theorem]{Example}

\theoremstyle{remark}
\newtheorem{remark}[theorem]{Remark}
\numberwithin{equation}{section}

\begin{document}

\setcounter{page}{1}

\title[The structure of Twisted power partial isometries]{The structure of Twisted power partial isometries}

\author[Athul Augustine \MakeLowercase{and} P. Shankar]{Athul Augustine \MakeLowercase{and} P. Shankar}

\address{Athul Augustine, Department of Mathematics, Cochin University of Science And Technology,  Ernakulam, Kerala - 682022, India. }
\email{\textcolor[rgb]{0.00,0.00,0.84}{athulaugus@gmail.com, athulaugus@cusat.ac.in}}

\address{P. Shankar, Department of Mathematics, Cochin University of Science And Technology,  Ernakulam, Kerala - 682022, India.}
\email{\textcolor[rgb]{0.00,0.00,0.84}{shankarsupy@gmail.com, shankarsupy@cusat.ac.in}}

\subjclass[2020]{Primary 47A13, 47A15, 46L65, 47A67, 47A20, 46L05, 81S05.}

\keywords{Power partial isometry; Halmos and Wallen decompositions; Reducing subspace; irreducible representation; non-commutative torus; Hardy space over the unit polydisc.}


\begin{abstract}
Let $n>1$ and let $\{U_{ij}\}_{1\leq i<j\leq n}$ be $n\choose 2$ commuting unitaries on a Hilbert space $\mathcal{H}$. Suppose $U_{ji}:=U^*_{ij}$, $1\leq i<j\leq n$. An n-tuple of power partial isometries $(V_1,...,V_n)$ on Hilbert space $\mathcal{H}$ is called $\mathcal{U}_n$-twisted power partial isometry with respect to $\{U_{ij}\}_{i<j}$ (or simply $\mathcal{U}_n$-twisted power partial isometry if $\{U_{ij}\}_{i<j}$ is clear from the context) if $V_i^*V_j=U_{ij}V_jV^*_i, ~~ V_iV_j=U_{ji}V_jV_i ~~\text{and}~~ V_kU_{ij}=U_{ij}V_k~~(i,j,k=1,2,...,n,~\text{and}~i\neq j).$ We prove that each $\mathcal{U}_n$-twisted power partial isometry admits a Halmos and Wallen \cite{HW70} type orthogonal decomposition.
\end{abstract}
\maketitle

\section{Introduction}
The Wold-von Neumann theorem states that every isometry on a Hilbert space is either a shift, a unitary, or a direct sum of shift and unitary. An operator $V$ is a \textit{power partial isometry} if $V^n$ is a partial isometry for all $n\geq 0$. Halmos and Wallen \cite{HW70} proved a similar result for power partial isometries as of the Wold-von Neumann theorem. Their theorem states that every power partial isometry is a direct sum of a unitary operator, some unilateral (forward) shifts, some backward shifts, and some truncated shifts on finite-dimensional spaces.

Slocinski \cite{MS80} proved an analogous result of the Wold-von Neumann theorem for a pair of doubly commuting isometries. Sarkar \cite{JS14} extended the ideas of Slocinski on the Wold-type decomposition for a pair of doubly commuting isometries to the multivariable case $(n\geq 2)$. Burdak \cite{BZ07} and Catepill\'{a}n and  Szyma\'{n}ski \cite{XW96} proved an analogous result of Halmos and Wallen for pairs of star-commuting power partial isometries. Heuf, Raeburn and Tolich \cite{HRT15} proved a Halmos and Wallen type structure theorem for a finite family of star-commuting (doubly commuting) power partial isometries.

Jeu and Pinto \cite{MP20} proved that n-tuple of doubly non-commuting isometries admits an orthogonal decomposition similar to the Wold-von Neumann type theorem. Rakshit, Sarkar, and Suryawanshi \cite{NJM22} extended the results of Jeu and Pinto for n-tuple of  $\mathcal{U}_n$-twisted isometries. Ostrovskyi, Proskurin, and Yakymiv \cite{OPY22} proved that an irreducible family of twisted commuting power partial isometries admits an orthogonal decomposition analogous to Halmos and Wallen theorem.

In this paper, we prove an orthogonal decomposition theorem for family of $\mathcal{U}_n$-twisted power partial isometries and analogous to the Halmos-Wallen type theorem. The paper is organized as follows. In Section 2, we set up notations and stated the Halmos and Wallen decomposition theorem for power partial isometry. In Section 3, we define $\mathcal{U}_n$-twisted power partial isometries and establish some preliminary results. 
In Section 4, we discuss some examples of $\mathcal{U}_n$-twisted power partial isometries. In Section 5, we prove Halmos and Wallen type decomposition theorem for irreducible family of  $\mathcal{U}_n$-twisted power partial isometries. In Section 6, we demonstrate the structure theorem for n-tuple of $\mathcal{U}_n$-twisted power partial isometries similar to Halmos and Wallen type decomposition theorem without irreducibility assumption.

\section{Power partial isometries}
In this section, we recall the Halmos and Wallen \cite{HW70} decomposition theorem for power partial isometries and some elementary facts about partial isometries and power partial isometries.

Let $V$ be a partial isometry on a Hilbert space $\mathcal{H}$. Suppose $\mathcal{K}$ is a subspace of $\mathcal{H}$ and $\mathcal{K}$ is reducing for $V$. Then $V|_{\mathcal{K}}$ is a partial isometry \cite[Lemma 1]{HW70}. If $V$ and $W$ are partial isometries on a Hilbert space $\mathcal{H}$, then $VW$ is a partial isometry if and only if the initial projection  $V^*V$ and range projection $WW^*$ commutes \cite[Lemma 2]{HW70}. An operator $V$ is said to be \textit{power partial isometry} if $V^n$ is a partial isometry for all $n\geq 0$. If $V$ is a power partial isometry then the family of projections $\{V^nV^{*n}\}\cup \{V^{*n}V^n\}$ will commute. (Notational convention: $V^{*n}=(V^*)^n$ and $V^{*n-m}=(V^*)^{n-m}$ with $n\geq m$.)

Unitary operators, the unilateral shift $S$ on $\ell^2$, the backward shift $S^*$ on $\ell^2$, and the truncated shifts $J_p$ on $\mathbb{C}^p$ are examples of power partial isometries. $J_p$ is defined as follows on the standard basis on $\mathbb{C}^p$, $J_p(e_n)=e_{n+1}$ for $n<p$ and $J_p(e_p)=0$. Note that $p\geq 1$ and $J_1=0$. Halmos and Wallen proved that every power partial isometry is a direct sum of these examples. 

\begin{theorem}\label{halwal}
(Halmos and Wallen). Let $V$ be a power partial isometry on a Hilbert space $\mathcal{H}$, and let $P$ and $Q$ be the orthogonal projections on $\cap_{n=1}^\infty V^n\mathcal{H}$ and $\cap_{n=1}^\infty V^{*n}\mathcal{H}$ respectively. Then $PQ=QP$ and the subspaces $\mathcal{H}_u:=PQ\mathcal{H}$, $\mathcal{H}_s:=(1-P)Q\mathcal{H}$, $\mathcal{H}_b:=(1-Q)P\mathcal{H}$ and 
$$\mathcal{H}_p:=\sum_{n=1}^p(V^{n-1}V^{*n-1}-V^nV^{*n})(V^{*p-n}V^{p-n}-V^{*p-n+1}V^{p-n+1})\mathcal{H} $$
are all reducing for $V$, and satisfy $\mathcal{H}=\mathcal{H}_u\oplus \mathcal{H}_s\oplus \mathcal{H}_b\oplus (\bigoplus_{p=1}^\infty \mathcal{H}_p)$. Further there are Hilbert spaces $M_s$, $M_b$ and $\{M_p:p\geq 1\}$ (allowing $M_{*}=\{0\}$) such that
\begin{enumerate}
\item[(a)] $V|_{\mathcal{H}_u}$ is unitary;
\item[(b)] $V|_{\mathcal{H}_s}$ is unitarily equivalent to $S\otimes 1$ on $\ell^2(\mathbb{N})\otimes M_s$;
\item[(c)] $V|_{\mathcal{H}_b}$ is unitarily equivalent to $S^*\otimes 1$ on $\ell^2(\mathbb{N})\otimes M_b$;
\item[(d)] for $p\geq 1$, $V|_{\mathcal{H}_p}$ is unitarily equivalent to $J_p\otimes 1$ on $\mathbb{C}^p\otimes M_p$.
\end{enumerate}
\end{theorem}

The \textit{multiplicity spaces $M_*$} are unique up to isomorphism, and thus the dimension of the multiplicity space is the only invariant. According to \cite{HRT15}, it is convenient to take multiplicity space to the subspaces mentioned as follows: \linebreak
$M_s=(1-VV^*)Q\mathcal{H}$, $M_b=(1-V^*V)P\mathcal{H}$, $M_p=(1-VV^*)(V^{*p-1}V^{p-1}-V^{*p}V^{p})\mathcal{H}$ for $p\geq 2$, and   $M_1=(1-VV^*)(1-V^*V)\mathcal{H}=\ker(V)\cap \ker(V^*)$.

From \cite{HW70, HRT15}, observe that the projections $V^nV^{*n}$ onto the subspaces $V^n\mathcal{H}$ form a decreasing sequence. Thus, $V^nV^{*n}$ converge to the projection $P$ onto $\cap_{n=1}^\infty V^n\mathcal{H}$ in the strong-operator topology \cite[Corollary 2.5.7]{KR97}. Also, projections $V^{*n}V^{n}$ onto the subspaces $V^{*n}\mathcal{H}$ form a decreasing sequence. Thus, $V^{*n}V^{n}$ converge to the projection $Q$ onto $\cap_{n=1}^\infty V^{*n}\mathcal{H}$ in the strong-operator topology. On the norm-bounded sets, the composition is jointly strong-operator continuous \cite[Remark 2.5.10]{KR97}. Thus $(V^nV^{*n})(V^{*n}V^n)$ converge to the projection $PQ$ onto\linebreak $\cap_{n=1}^\infty V^n\mathcal{H} \cap \cap_{n=1}^\infty V^{*n}\mathcal{H}$ in the strong-operator topology. All the range and source projections commute, thus $PQ=QP$. Therefore all the product of projections $PQ, ~(1-P)Q$, etc., are projections onto the respective subspaces. Also, the subspaces of $\mathcal{H}$ corresponding to each projection are closed.

\section{$\mathcal{U}_n$-twisted power partial isometries}

Let $\lambda_{ij}\in \mathbb{T}, 1\leq i < j\leq n$, and suppose that $\lambda_{ji}=\overline{\lambda_{ij}}$ for all $1\leq i < j\leq n$. A family of power partial isometries $(V_1,....,V_n)$ on some Hilbert space $\mathcal{H}$ is said to be \textit{twisted commuting power partial isometries} if ${V_i}^*V_j=\lambda_{ij}V_j{V_i}^*$ and $V_iV_j=\lambda_{ji}V_jV_i$ for all $i\neq j$. Each irreducible family of twisted commuting power partial isometries admits an orthogonal decomposition analogous to Halmos and Wallen theorem \cite{OPY22}.

If $\lambda_{ij}=1, i\neq j$, then the twisted commuting power partial isometries are simply doubly commuting isometries. Then the condition reduces to orthogonal decompositions of star-commuting power partial isometries \cite{HRT15}. A question of apparent  interest is to enlarge the above class of family of power partial isometries that admit the orthogonal decomposition. To answer this question, we now introduce our primary object of study, $\mathcal{U}_n$-twisted power partial isometries on Hilbert spaces.

\begin{definition}
\textbf{($\mathcal{U}_n$-twisted power partial  isometries).}  Let $n>1$ and let $\{U_{ij}\}_{1\leq i<j\leq n}$ be $n\choose 2$ commuting unitaries on a Hilbert space $\mathcal{H}$. Suppose $U_{ji}:=U^*_{ij}$, $1\leq i<j\leq n$. A family of power partial isometries $(V_1,...,V_n)$ on Hilbert space $\mathcal{H}$ is called $\mathcal{U}_n$-twisted power partial isometry with respect to $\{U_{ij}\}_{i<j}$ if 
\begin{equation*}
V_i^*V_j=U_{ij}V_jV^*_i, ~~ V_iV_j=U_{ji}V_jV_i ~~\text{and}~~ V_kU_{ij}=U_{ij}V_k~~(i,j,k=1,2,...,n,~\text{and}~i\neq j).
\end{equation*} 
\end{definition}

Sometimes we will simply say that $(V_1,...,V_n)$ is a $\mathcal{U}_n$-twisted power partial isometry without referring the unitaries $\{U_{ij}\}_{1\leq i<j\leq n}$. Two $\mathcal{U}_n$-twisted power partial  isometries $(V_1,...,V_n)$ on the Hilbert space $\mathcal{H}$ and $(W_1,...,W_n)$ on the Hilbert space $\mathcal{K}$ are said to be simultaneously unitarily equivalent if there is a unitary isomorphism $U$ from $\mathcal{H}$ onto $\mathcal{K}$ such that $UV_iU^*=W_i$ for all $1\leq i\leq n$. Clearly, twisted commuting power partial isometries are also $\mathcal{U}_n$-twisted power partial isometries with respect to $\{\lambda_{ij}I_\mathcal{H}\}_{i<j}$.

\begin{lemma}\label{U2TWIST}
Let $(V,W)$ be a $\mathcal{U}_2$-twisted power partial  isometries on a Hilbert space $\mathcal{H}$. Let $P$ and $Q$ be the projections on the subspaces $\cap_{n=1}^\infty V^n\mathcal{H}$ and $\cap_{n=1}^\infty V^{*n}\mathcal{H}$ respectively. Then $P$ and $Q$ are $\mathcal{U}_2$-twisted with $W$.
\end{lemma}
\begin{proof}
Since $(V,W)$ is $\mathcal{U}_2$-twisted power partial  isometries  then there exists a unitary $U \in B(\mathcal{H})$ such that $V^*W=UWV^*$, $VW=U^*WV$ and $V,W\in \{U\}'$. For all $n\geq 1$, we have 
\begin{equation*}
\begin{split}
(V^nV^{*n})W &= V^nV^{*n-1}UWV^*\\
&= V^nV^{*n-2}UV^*WV^*\\
&= V^nU^{n}WV^{*n}\\
&= U^{n}V^nWV^{*n}\\
&= U^{n}U^{n*}W(V^nV^{*n})^*\\
&= IW(V^nV^{*n})^*.\\
\end{split}
\end{equation*}
and
\begin{equation*}
\begin{split}
(V^{*n}V^{n})W &= V^{*n}V^{n-1}U^*WV\\
&= V^{*n}U^{*n}WV^{*n}\\
&= U^{*n}V^nWV^{*n}\\
&= U^{*n}U^{n}W(V^nV^{*n})^*\\
&= IW(V^nV^{*n})^*.\\
\end{split}
\end{equation*}
Since the projections $P$ and $Q$ are strong operator limits of the sequences $\{V^nV^{*n}\}$ and $\{V^{*n}V^n\}$, respectively. It follows that $P$ and $Q$ are $\mathcal{U}_2$-twisted with $W$.
\end{proof}

\begin{lemma}\label{PQtwist}
Let $(V,W)$ be a $\mathcal{U}_2$-twisted power partial  isometries on a Hilbert space $\mathcal{H}$. Let $P$ and $Q$ be the projections on the subspaces $\cap_{n=1}^\infty V^n\mathcal{H}$ and $\cap_{n=1}^\infty V^{*n}\mathcal{H}$ respectively. Then the projections $PQ,~(1-P)Q$ and $(1-Q)P$ are $\mathcal{U}_2$-twisted with $W$.
\end{lemma}
\begin{proof}
From Lemma \ref{U2TWIST}, the projections $P$ and $Q$ are $\mathcal{U}_2$-twisted with $W$. We have
\begin{equation*}
\begin{split}
Q(1-P)W &=QW-QPW\\
        &=IWQ-IWQP\\
        &=WQ(1-P)\\
        &=WQ(1-P).\\
\end{split}
\end{equation*}
and 
\begin{equation*}
\begin{split}
(1-P)QW &=QW-PQW\\
        &=IWQ-IWPQ\\
        &=W(1-P)Q\\
        &=W(1-P)Q.\\
\end{split}
\end{equation*}
Thus projection $(1-P)Q$ is $\mathcal{U}_2$-twisted with $W$. Similarly, we can prove that the projections $(1-Q)P$ and $PQ$ are $\mathcal{U}_2$-twisted with $W$.
\end{proof}

\begin{lemma}\label{Hptwist}
Let $(V,W)$ be a $\mathcal{U}_2$-twisted power partial  isometries on a Hilbert space $\mathcal{H}$. Let 
$\mathcal{H}_p=\sum_{n=1}^p(V^{n-1}V^{*n-1}-V^nV^{*n})(V^{*p-n}V^{p-n}-V^{*p-n+1}V^{p-n+1})\mathcal{H} $ be a subspace of $\mathcal{H}$. Then the projection onto $\mathcal{H}_p$ is $U_2$-twisted with $W$. 
\end{lemma}
\begin{proof}
Since projection onto $\mathcal{H}_p$ involves only range and source projections of $V^n$. By Lemma \ref{U2TWIST}, the source and the range projections of $V^n$ are $\mathcal{U}_2$-twisted with $W$. Hence the projection onto $\mathcal{H}_p$ is $\mathcal{U}_2$-twisted with $W$.
\end{proof}

\begin{lemma}\label{reducing}
Let $\{V_1,....,V_n\}$ be a family of $\mathcal{U}_n$-twisted power partial  isometries on a Hilbert space $\mathcal{H}$. Let $\mathcal{H}$ has a decomposition of $V_1$,
$$\mathcal{H}=\mathcal{H}_u\oplus \mathcal{H}_s\oplus \mathcal{H}_b\oplus (\bigoplus_{p=1}^\infty \mathcal{H}_p)$$
as in Theorem \ref{halwal}, where the subspaces $\mathcal{H}_u,~\mathcal{H}_s,~\mathcal{H}_b$ and $\mathcal{H}_p,~p\geq 1$ are reducing for $V_1$. Then the subspaces $\mathcal{H}_u,~\mathcal{H}_s,~\mathcal{H}_b$ and $\mathcal{H}_p,~p\geq 1$ are reducing for $V_i$, $i=2,...,n.$
\end{lemma}

\begin{proof}
Let $P$ and $Q$ be the orthogonal projections onto $\cap_{n=1}^\infty V_1^n\mathcal{H}$ and $\cap_{n=1}^\infty V_1^{*n}\mathcal{H}$ respectively.
By Thereom \ref{halwal}, we have $\mathcal{H}_u= PQ\mathcal{H}$. For $k\in \mathcal{H}_u$, $k=PQh$ for some $h\in \mathcal{H}$. From Lemma \ref{PQtwist}, $PQ$ is $\mathcal{U}_n$-twisted with $V_i$, then
$$V_ik=V_iPQh=UPQV_ih=PQUV_ih\in PQ\mathcal{H}=\mathcal{H}_u $$
and 
$$V_i^*k=V_i^*PQh=UPQV_i^*h=PQUV_i^*h\in PQ\mathcal{H}=\mathcal{H}_u.$$
Thus $\mathcal{H}_u$ is reducing for $V_i$. Similarly, we can prove $\mathcal{H}_s$, $\mathcal{H}_b$ and $\mathcal{H}_p,~p\geq 1$ are reducing for $V_i$.
\end{proof}

\begin{corollary}
Let $\{V_1,....,V_n\}$ be an irreducible family of $\mathcal{U}_n$-twisted power partial  isometries on a Hilbert space $\mathcal{H}$. Then $\mathcal{H}$ coincides with exactly one of the components of its Halmos and Wallen decomposition theorem.
\end{corollary}

\section{Examples}
In this section,  we discuss some basic concepts and present some (model) examples of $\mathcal{U}_n$-twisted power partial isometries. This section takes a comprehensive approach to $\mathcal{U}_n$-twisted power partial isometries in what follows. This section is the core part of this paper. The examples are motivated by the ideas from Rakshit, Sarkar, and Suryawanshi in \cite{NJM22}.

Let $H^2(\mathbb{D})$ denote the Hardy space over the unit disc\linebreak $\mathbb{D}=\{z\in \mathbb{C} : |z|<1\}.$ Then the multiplication operator on $H^2(\mathbb{D})$ by the coordinate function $z$ is denoted by $M_zf=zf$ for all $f\in H^2(\mathbb{D})$. It is easy to observe that $M_z$ is a shift operator on $H^2(\mathbb{D})$ of multiplicity one (as ker$M_z^*=\mathbb{C}$). Now let $H^2(\mathbb{D}^2)$ be the Hardy space over the bidisc $\mathbb{D}^2.$ Then $H^2(\mathbb{D}^2)$ is the Hilbert space of all square summable analytic functions on $\mathbb{D}^2.$ An analytic function $f(z) = \sum_{k\in \mathbb{Z}_+^2} \alpha_kz^k$  on $\mathbb{D}^2$ is in $H^2(\mathbb{D}^2)$ if and only if 
$$ ||f||:= \left( \sum_{k\in \mathbb{Z}_+^2} |\alpha_k|^2\right)^\frac{1}{2} < \infty.  $$
One can easily identify $H^2(\mathbb{D}^2)$ with $H^2(\mathbb{D})\otimes H^2(\mathbb{D})$ in a natural way by defining $\sigma : H^2(\mathbb{D})\otimes H^2(\mathbb{D}) \rightarrow H^2(\mathbb{D}^2)$  by $\sigma (z^{k_1}\otimes z^{k_2}) = z_{1}^{k_1} z_{2}^{k_2},~ k\in \mathbb{Z}_+^2.$ Then $\sigma$ is a unitary operator and
$$ \sigma(M_z\otimes I_{H^2(\mathbb{D})}) = M_{z_1}\sigma \qquad\text{and}\qquad\sigma(I_{H^2(\mathbb{D})}\otimes M_z) = M_{z_2}\sigma,$$
where $M_{z_1}$ and $M_{z_2}$ are the multiplication operators by $z_1$ and $z_2$, respectively, on $H^2(\mathbb{D}^2)$. The above construction of $H^2(\mathbb{D}^2)$ works equally well for the Hardy space $H^2(\mathbb{D}^m)$ over $\mathbb{D}^m,~ m>1$.

\begin{example}
We now introduce a special class of diagonal operators parametrized by the circle group $\mathbb{T}$. For each $\lambda\in \mathbb{T}$, define 
$$ D[\lambda]z^m = \lambda^mz^m \qquad (m\in \mathbb{Z}_+).$$
$D[\lambda]$ is a unitary diagonal operator on $H^2(\mathbb{D})$ and $D[\lambda]^* = D[\overline{\lambda}] = $ diag$(1,\overline{\lambda},\overline{\lambda}^2,....)$. It is easy to observe that,
\[
    (M_z^*D[\lambda])(z^m)= \begin{cases}
  \lambda^mz^{m-1}  &  \text{if}~~ m>0, \\
  0 & \text{if}~~ m=0,
\end{cases}
  \]
and 
\[
    (D[\lambda]M_z^*)(z^m)= \begin{cases}
  \lambda^{m-1}z^{m-1}  &  \text{if}~~ m>0, \\
  0 & \text{if}~~ m=0.
\end{cases}
  \]
Also
$$ (M_zD[\lambda])(z^m)= \lambda^m z^{m+1} \qquad\text{for}\qquad m\in \mathbb{Z}_+ $$ 
\text{and} 
$$  (D[\lambda]M_z)(z^m)= \lambda^{m+1}z^{m+1} \qquad \text{for}\qquad  m\in \mathbb{Z}_+.$$
Therefore, we have $M_z^*D[\lambda] = \lambda D[\lambda]M_z^*$ and $D[\lambda]M_z = \lambda M_zD[\lambda]$. Now, we fix $\lambda\in \mathbb{T}$, and define $T_1$ and $T_2$ on $H^2(\mathbb{D}^2)$ as
$$ T_1 = M_z\otimes I_{H^2(\mathbb{D})}\qquad \text{and} \qquad T_2 = D[\lambda] \otimes M_z.$$
Then, $(T_1,T_2)$ is a pair of power partial isometries on $H^2(\mathbb{D}^2)$. We have\linebreak  $T_1^*T_2 = M_z^*D[\lambda]\otimes M_z$ and $T_2T_1^*= D[\lambda]M_z^*\otimes M_z$. $M_z^*D[\lambda] = \lambda D[\lambda]M_z^*$ implies $T_1^*T_2=\lambda T_2T_1^*$. Also, $T_1T_2 = M_zD[\lambda]\otimes M_z$ and $T_2T_1 = D[\lambda]M_z \otimes M_z$. Then, $D[\lambda]M_z = \lambda M_zD[\lambda]$ implies that $T_2T_1 = \lambda T_1T_2$.

We now consider the Hilbert space $\mathcal{H} = H^2(\mathbb{D}^2) \oplus H^2(\mathbb{D}^2)$, and the power partial isometries $V_1 =$ diag$(T_1,T_2)$ and $V_2 =$diag$(T_2,T_1)$ on $\mathcal{H}$. If we set \linebreak $U =$ diag$(\lambda I_{H^2(\mathbb{D}^2)},\overline{\lambda} I_{H^2(\mathbb{D}^2)})$, then

$$V_1^*V_2 = \begin{bmatrix}
T_1^*T_2 & 0 \\
0 & T_2^*T_1 
\end{bmatrix}
= \begin{bmatrix}
\lambda T_2T_1^* & 0 \\
0 & \overline{\lambda}T_1T_2^* 
\end{bmatrix}
= \begin{bmatrix}
\lambda I_{H^2(\mathbb{D}^2)} & 0 \\
0 & \overline{\lambda} I_{H^2(\mathbb{D}^2)} 
\end{bmatrix}V_2V_1^*$$
 and 
$$V_1V_2 = \begin{bmatrix}
T_1T_2 & 0 \\
0 & T_2T_1 
\end{bmatrix}
= \begin{bmatrix}
\overline{\lambda} T_2T_1 & 0 \\
0 & \lambda T_1T_2 
\end{bmatrix}
= \begin{bmatrix}
\overline{\lambda} I_{H^2(\mathbb{D}^2)} & 0 \\
0 & \lambda I_{H^2(\mathbb{D}^2)} 
\end{bmatrix}V_2V_1$$ 
which implies that $V_1^*V_2 = UV_2V_1^*$ and $V_1V_2 = U^*V_2V_1$. Since $V_1,V_2 \in \{U\}^{'}$, thus  the pair $(V_1,V_2)$ is a $\mathcal{U}_2$-twisted power partial isometry on $\mathcal{H}$.
\end{example}

\begin{example}
Let $\lambda\in \mathbb{T}$ and $D[\lambda]$ as in above example. Define $T_3$ and $T_4$ in $H^2(\mathbb{D}^2)$ as
$$ T_3 = M_z^*\otimes I_{H^2(\mathbb{D})}\qquad \text{and} \qquad T_4 = D[\lambda] \otimes M_z^*.$$
Then, $(T_3,T_4)$ is a pair of power partial isometries on $H^2(\mathbb{D}^2)$. We have \linebreak $T_3^*T_4 = M_zD[\lambda]\otimes M_z^*$ and $T_4T_3^*= D[\lambda]M_z\otimes M_z^*$. $D[\lambda]M_z = \lambda M_zD[\lambda]$ implies $T_4T_3^*= \lambda T_3^*T_4$. Also, $T_3T_4 = M_z^*D[\lambda]\otimes M_z^*$ and $T_4T_3 = D[\lambda]M_z^* \otimes M_z^*$. Then, $M_z^*D[\lambda] = \lambda D[\lambda]M_z^*$ implies that $T_3T_4 = \lambda T_4T_3$.

We now consider the Hilbert space $\mathcal{H} = H^2(\mathbb{D}^2) \oplus H^2(\mathbb{D}^2)$, and the power partial isometries $V_3 =$ diag$(T_3,T_4)$ and $V_4 =$diag$(T_4,T_3)$ on $\mathcal{H}$. If we set \linebreak $U =$ diag$(\overline{\lambda} I_{H^2(\mathbb{D}^2)},\lambda I_{H^2(\mathbb{D}^2)})$, then

$$V_3^*V_4 = \begin{bmatrix}
T_3^*T_4 & 0 \\
0 & T_4^*T_3 
\end{bmatrix}
= \begin{bmatrix}
\overline{\lambda} T_4T_3^* & 0 \\
0 & \lambda T_3T_4^* 
\end{bmatrix}
= \begin{bmatrix}
\overline{\lambda} I_{H^2(\mathbb{D}^2)} & 0 \\
0 & \lambda I_{H^2(\mathbb{D}^2)} 
\end{bmatrix}V_4V_3^*$$
 and 
$$V_3V_4 = \begin{bmatrix}
T_3T_4 & 0 \\
0 & T_4T_3 
\end{bmatrix}
= \begin{bmatrix}
\lambda T_4T_3 & 0 \\
0 & \overline{\lambda} T_3T_4 
\end{bmatrix}
= \begin{bmatrix}
\lambda I_{H^2(\mathbb{D}^2)} & 0 \\
0 & \overline{\lambda} I_{H^2(\mathbb{D}^2)} 
\end{bmatrix}V_4V_3$$ 
which implies that $V_3^*V_4 = UV_4V_3^*$ and $V_3V_4 = U^*V_4V_3$. Since $V_3,V_4 \in \{U\}^{'}$, thus the pair $(V_3,V_4)$ is a $\mathcal{U}_2$-twisted power partial isometry on $\mathcal{H}$.
\end{example}

\begin{example}
Let $J_p$ be the truncated shifts on $\mathbb{C}^p$ and define  $J_p$ in terms of the standard basis on $\mathbb{C}^p$  by $J_pe_n = e_{n+1}$ for $n < p$ and $J_pe_p=0$. Note that $p\geq 1$ and $J_1=0$. For each $\lambda\in \mathbb{T}$, define

$$
   d[\lambda]e_m = \lambda^{m-1}e_m.
  $$
Clearly,  $d[\lambda]$ is a unitary diagonal operator on $\mathbb{C}^p$ and 
$$d[\lambda]^* = d[\overline{\lambda}] =  \text{diag}(1,\overline{\lambda},\overline{\lambda}^2,...,\overline{\lambda}^{p-1}).$$ It is easy to see that 
\[
    (J_p^*d[\lambda])(e_m)= \begin{cases}
  \lambda^{m-1}e_{m-1}  &  \text{if}~~ m>1, \\
  0 & \text{if}~~ m=1,
\end{cases}
  \]
and 
\[
    (d[\lambda]J_p^*)(e_m)= \begin{cases}
 \lambda^{m-2}e_{m-1}  &  \text{if}~~ m>1, \\
  0 & \text{if}~~m=1.
\end{cases}
  \]
Also
\[
    (J_pd[\lambda])(e_m)= \begin{cases}
  \lambda^{m-1}e_{m+1}  &  \text{if}~~  m<p, \\
  0 & \text{if}~~ m=p,
\end{cases}
  \]
and 
\[
    (d[\lambda]J_p)(e_m)= \begin{cases}
  \lambda^{m}e_{m+1}  &  \text{if}~~  m<p, \\
  0 & \text{if}~~ m=p.
\end{cases}
  \]
Thus, we have $J_p^*d[\lambda] = \lambda d[\lambda]J_p^*$ and $d[\lambda]J_p = \lambda J_pd[\lambda]$. Now fix $\lambda\in \mathbb{T}$, and define $S_1$ and $S_2$ in $\mathbb{C}^p\otimes \mathbb{C}^p$ as
$$ S_1 = J_p\otimes I_{\mathbb{C}^p}\qquad \text{and} \qquad S_2 = d[\lambda] \otimes J_p.$$
Then, $(S_1,S_2)$ is a pair of power partial isometries on $\mathbb{C}^p\otimes \mathbb{C}^p$, and \linebreak $S_1^*S_2 = J_p^*d[\lambda]\otimes J_p$ and $S_2S_1^*= d[\lambda]J_p^*\otimes J_p$. $J_p^*d[\lambda] = \lambda d[\lambda]J_p^*$ implies have $S_1^*S_2=\lambda S_2S_1^*$. Also, $S_1S_2 = J_pd[\lambda]\otimes J_p$ and $S_2S_1 = d[\lambda]J_p \otimes J_p$. Then, $d[\lambda]J_p = \lambda J_pd[\lambda]$ implies that $S_2S_1 = \lambda S_1S_2$.

We now consider the Hilbert space $\mathcal{H} = (\mathbb{C}^p\otimes \mathbb{C}^p) \oplus (\mathbb{C}^p\otimes \mathbb{C}^p)$, and the power partial isometries $M_1 =$ diag$(S_1,S_2)$ and $M_2 =$diag$(S_2,S_1)$ on $\mathcal{H}$. If we set \linebreak $U =$ diag$(\lambda I_{H^2(\mathbb{D}^2)},\overline{\lambda} I_{H^2(\mathbb{D}^2)})$, then

$$M_1^*M_2 = \begin{bmatrix}
S_1^*S_2 & 0 \\
0 & S_2^*S_1 
\end{bmatrix}
= \begin{bmatrix}
\lambda S_2S_1^* & 0 \\
0 & \overline{\lambda}S_1S_2^* 
\end{bmatrix}
= \begin{bmatrix}
\lambda I_{H^2(\mathbb{D}^2)} & 0 \\
0 & \overline{\lambda} I_{H^2(\mathbb{D}^2)} 
\end{bmatrix}M_2M_1^*$$
 and 
$$M_1M_2 = \begin{bmatrix}
S_1S_2 & 0 \\
0 & S_2S_1 
\end{bmatrix}
= \begin{bmatrix}
\overline{\lambda} S_2S_1 & 0 \\
0 & \lambda S_1S_2 
\end{bmatrix}
= \begin{bmatrix}
\overline{\lambda} I_{H^2(\mathbb{D}^2)} & 0 \\
0 & \lambda I_{H^2(\mathbb{D}^2)} 
\end{bmatrix}M_2M_1$$ 
which implies that $M_1^*M_2 = UM_2M_1^*$ and $M_1M_2 = U^*M_2M_1$. Since \linebreak $M_1,M_2 \in \{U\}^{'}$, thus the pair $(M_1,M_2)$ is a $\mathcal{U}_2$-twisted power partial isometry on $\mathcal{H}$.
\end{example}

It is clear that for each $\lambda\in \mathbb{T}$, the pairs $(M_z,D[\lambda])$, $(T_1,T_2)$, $(M_z^*,D[\lambda])$, $(T_3,T_4)$, $(J_p,d[\lambda])$ and $(S_1,S_2)$, defined in above examples are twisted commuting power partial isometries \cite{OPY22}.

We now extend the discussion of Hardy space over $\mathbb{D}^m, m>1$. Let $\mathcal{E}$ be a Hilbert space. Let $H_{\mathcal{E}}^2(\mathbb{D}^m)$ denote the $\mathcal{E}$-valued Hardy space over $\mathbb{D}^m$. Observe that $H_{\mathcal{E}}^2(\mathbb{D}^m)$ is the Hilbert space of all square summable analytic functions on $\mathbb{D}^m$ with coefficients in $\mathcal{E}$. So, if $\mathcal{E}=\mathbb{C}$ then $H^2(\mathbb{D}^m)=H_{\mathbb{C}}^2(\mathbb{D}^m)$. By natural identification 
$$ z^k\eta \leftrightarrow z^{k_1}\otimes ....\otimes z^{k_m}\otimes \eta \qquad (k\in \mathbb{Z}_+^m, \eta\in \mathcal{E}),$$
up to unitary equivalence, we have
$$H_{\mathcal{E}}^2(\mathbb{D}^m) = \underbrace{H^2(\mathbb{D})\otimes ....\otimes H^2(\mathbb{D})}_\textrm{$m$-times} \otimes \mathcal{E} = H^2(\mathbb{D}^m)\otimes \mathcal{E}.$$
In the above setting, for each fixed $i=1,2,....,m$, up to unitary equivalence, we have 
$$M_{z_i}= (I_{H^2(\mathbb{D})}\otimes ...I_{H^2(\mathbb{D})}\otimes \underbrace{M_z}_\textrm{$i$-th}  \otimes I_{H^2(\mathbb{D})}\otimes...\otimes I_{H^2(\mathbb{D})})\otimes I_{\mathcal{E}} = M_{z_i} \otimes I_{\mathcal{E}},$$
where $M_{z_i}f=z_if$ for any $f$ in $H_{\mathcal{E}}^2(\mathbb{D}^m)$. We shall use the above identification interchangeably whenever appropriate. The above tensor product representation of the multiplication operators implies that $(M_{z_1},M_{z_2},...M_{z_m})$ on $H_{\mathcal{E}}^2(\mathbb{D}^m)$ is \textit{doubly commuting}, that is, $M_{z_i}M_{z_j}=M_{z_j}M_{z_i}$ and $M_{z_i}^*M_{z_p}=M_{z_p}M_{z_i}^*$ for all $i,j,p=1,2,...,m$ and $i\neq p$.

Also, consider 
$$ \underbrace{\mathbb{C}^p \otimes .....\otimes \mathbb{C}^p}_\textrm{m times} \otimes \mathcal{E} $$
and analogue to above setting, for each fixed $i= 1,2,...,m$, up to unitary equivalence, we have  
$$ J_{p_i} = (I_{\mathbb{C}^p} \otimes...I_{\mathbb{C}^p}\otimes \underbrace{J_p}_\textrm{$i$-th}  \otimes I_{\mathbb{C}^p}\otimes...\otimes I_{\mathbb{C}^p})\otimes I_{\mathcal{E}} = J_{p_i} \otimes I_{\mathcal{E}}. $$

We present the vital notion $j$th diagonal operator defined by Rakshit, Sarkar, and Suryawanshi \cite{NJM22}. 

\begin{definition}
Let $\mathcal{E}$ be a given Hilbert space and $U\in B(\mathcal{E})$ be a unitary operator. For $j\in \{1,...,m\}$, the \textit{$jth$ diagonal operator} $D_j[U]$ with symbol $U$ on $H_{\mathcal{E}}^2(\mathbb{D}^m)$ is defined by
$$D_j[U](z^k\eta) = z^k(U^{k_j}\eta) \qquad(k\in \mathbb{Z}_+^m, \eta \in \mathcal{E}).$$
Here $k=(k_1,k_2,...,k_m).$ Note that $D_j[U]$ is a unitary operator on $H_{\mathcal{E}}^2(\mathbb{D}^m)$. In particular, if $m=1$ and $\mathcal{E}=\mathbb{C}$, then $U$ is given by $U=\lambda$ for some $\lambda\in \mathbb{T}$, and $D_1[\lambda]$ is the diagonal operator diag$(1,\lambda,\lambda^2,....)$ on $H^2(\mathbb{D})$.
\end{definition}

The following Lemma proved in \cite[Lemma 2.3]{NJM22} is helpful for us to study the decomposition of $\mathcal{U}_n$-twisted power partial isometries.

\begin{lemma}\label{mz}
Let $\mathcal{E}$ be a Hilbert space, and let $U$ and $\tilde{U}$ be commuting unitaries in $\mathcal{B}(\mathcal{E})$. Suppose $i,j \in \{1,...,n\}$. Then
\begin{enumerate}
\item $D_j[U]^* = D_j[U^*]$ and $D_i[U]D_j[\tilde{U}] = D_j[\tilde{U}]D_i[U].$
\item $M_{z_i} D_j[U] = D_j[U]M_{z_i}$ whenever $i\neq j$.
\item $M_{z_i}^* D_i[U] = (I_{H^2(\mathbb{D}^n)} \otimes U) D_i[U]M_{z_i}^*$.
\end{enumerate}
\end{lemma}

Now, we define the $j$th diagonal operator on $\mathbb{C}^p\otimes...\otimes \mathbb{C}^p\otimes \mathcal{E}$ (m times of $\mathbb{C}^p$) to proceed. 

\begin{definition}
Let $j\in \{1,...,m\}$ and $A =\{1,2,...,p\}$. Given a Hilbert space $\mathcal{E}$ and a unitary $U \in \mathcal{B}(\mathcal{E})$, the $jth$ diagonal operator with symbol $U$ is the unitary operator $d_j[U]$ on $\mathbb{C}^p\otimes...\otimes \mathbb{C}^p\otimes \mathcal{E}$ (m times of $\mathbb{C}^p$) defined by
$$d_j[U](e_k\eta) = e_k(U^{k_j-1}\eta) \qquad(k\in A^m, \eta \in \mathcal{E}).$$
Here $k=(k_1,k_2,...,k_m).$ In particular, if $m=1$ and $\mathcal{E}=\mathbb{C}$, then $U$ is given by $U=\lambda$ for some $\lambda\in \mathbb{T}$, and $d_1[\lambda]$ is the diagonal operator diag$(1,\lambda,\lambda^2,...,\lambda^{p-1})$ on $\mathbb{C}^p$.
\end{definition}

\begin{lemma}\label{jpd}
Let $\mathcal{E}$ be a Hilbert space. Let $U$ and $\tilde{U}$ be commuting unitaries in $\mathcal{B}(\mathcal{E})$. Suppose $i,j \in \{1,...,n\}$. Then
\begin{enumerate}
\item $d_j[U]^* = d_j[U^*]$ and $d_i[U]d_j[\tilde{U}] = d_j[\tilde{U}]d_i[U].$
\item $J_{p_i} d_j[U] = d_j[U]J_{p_i}$ whenever $i\neq j$.
\item $J_{p_i}^* d_i[U] = (I_{(\mathbb{C}^p\otimes ...\otimes \mathbb{C}^p)} \otimes U) d_i[U]J_{p_i}^*$.
\end{enumerate}
\end{lemma}
\begin{proof}
The first part follows from the definition of $d_j[U]$ and the commutativity of $U$ and $\tilde{U}$. To prove $(2)$, assume that $k\in A^m$ and let $\eta\in \mathcal{E}$. Suppose $i\neq j$. Then we have 

\[
    (J_{p_i} d_j[U])(e_k\eta) = J_{p_i}e_k(U^{k_j-1}\eta) =\begin{cases}
  e_{k+e_i}(U^{k_j-1}\eta)  &  \text{if}~~ k_i<p, \\
  0 & \text{if}~~ k_i=p,
\end{cases}
  \]
  and
 \[
    (d_j[U]J_{p_i})(e_k\eta) =\begin{cases}
  d_j[U](e_{k+e_i}\eta)  &  \text{if}~~ k_i<p, \\
  0 & \text{if}~~ k_i=p,
\end{cases} =\begin{cases}
  e_{k+e_i}(U^{k_j-1}\eta)  &  \text{if}~~ k_i<p, \\
  0 & \text{if}~~ k_i=p,
\end{cases}
  \] 
  where $e_i$ denotes the element in $A^m$ with $1$ in the $ith$ slot and zero elsewhere. The condition $i\neq j$ implies that $k_j$ remains unchanged.
  
  For the third part, using $d_i[U](e_{k-e_i}\eta)= e_{k-e_i}(U^{k_i-2}\eta)$ for $k_i > 1$, we compute
  $$(J_{p_i}^* d_i[U])(e_k\eta) = J_{p_i}^*e_k(U^{k_i-1}\eta) = \begin{cases}
  e_{k-e_i}(U^{k_i-1}\eta)  &  \text{if}~~ k_i>1, \\
  0 & \text{if}~~ k_i=1,
\end{cases}$$ 
and 
$$(d_i[U]J_{p_i}^*)(e_k\eta) = \begin{cases}
  d_i[U](e_{k-e_i}\eta)  &  \text{if}~~ k_i>1, \\
  0 & \text{if}~~ k_i=1,
\end{cases}= \begin{cases}
   e_{k-e_i}(U^{k_i-2}\eta) &  \text{if}~~ k_i>1, \\
  0 & \text{if}~~ k_i=1.
\end{cases}$$
 This completes the proof of part $(3).$
\end{proof}

Now we provide more general examples of $\mathcal{U}_n$-twisted power partial isometries.

\begin{proposition}
Let $\mathcal{E}$ be a Hilbert space, and let $\{U_{ij} : i,j = 1,2,...,n, i\neq j\}$ be a commuting family of unitaries on $\mathcal{E}$ such that $U_{ji}:=U_{ij}^*$ for all $i\neq j$. Fix $k,l,m \in \{1,....,n\}$ with $k<l<m$ and consider $(n-m)$ unitary operators $\{U_{m+1},....,U_n\}$ in $\mathcal{B}(\mathcal{E})$ such that 
$$ U_iU_j=U_{ij}U_jU_i\qquad\text{and}\qquad U_iU_{pq}=U_{pq}U_i$$
for all $m+1\leq i\neq j \leq n,$ and $1\leq p\neq q\leq n.$ Let $M_1=M_{z_1}$ and 

$$
   M_i =\begin{cases}
  \bigotimes_{j=1}^{i-1}D_j[U_{ij}]\otimes M_{z_i}\bigotimes_{i+1}^{k} 1_{H^2(\mathbb{D})}\bigotimes_{k+1}^{l} 1_{H^2(\mathbb{D})}\bigotimes_{l+1}^{m}1_{\mathbb{C}^p} \otimes 1_{\mathcal{E}},  &  \text{if}~~ 2\leq i \leq k, \\
  \bigotimes_{j=1}^{k}D_j[U_{ij}]\bigotimes_{j=k+1}^{i-1}D_j[U_{ji}]\otimes M_{z_i}^*\bigotimes_{i+1}^{l} 1_{H^2(\mathbb{D})}\bigotimes_{l+1}^{m}1_{\mathbb{C}^p} \otimes 1_{\mathcal{E}},  &  \text{if}~~ k+1\leq i\leq l, \\
  \bigotimes_{j=1}^{k}D_j[U_{ij}]\bigotimes_{j=k+1}^{l}D_j[U_{ji}]\bigotimes_{j=l+1}^{i-1}d_j[U_{ij}]\otimes J_{p_i}\bigotimes_{i+1}^{m}1_{\mathbb{C}^p} \otimes 1_{\mathcal{E}}, & \text{if}~~ l+1\leq i \leq m,\\
  \bigotimes_{j=1}^{k}D_j[U_{ij}]\bigotimes_{j=k+1}^{l}D_j[U_{ji}]\bigotimes_{j=l+1}^{m}d_j[U_{ij}]\otimes 1_{\mathcal{E}} \otimes U_i, & \text{if}~~ m+1\leq i \leq n.
\end{cases}
$$
Then $(M_1,...,M_k)$ are shifts, $(M_{k+1},...,M_l)$ are backward shifs, $(M_{l+1},...,M_m)$ are truncated shifts, $(M_{m+1},...,M_n)$ are unitaries, and $(M_1,...,M_n)$ is a $\mathcal{U}_n$-twisted power partial isometry on $\mathcal{H}$, where 
$$\mathcal{H}= H^2(\mathbb{D}^k) \otimes H^2(\mathbb{D}^{l-k}) \otimes \underbrace{\mathbb{C}^p \otimes .....\otimes \mathbb{C}^p}_\textrm{m-l times}\otimes \mathcal{E}$$
 with respect to $\{I_{\mathcal{K}} \otimes U_{ij}\}_{i<j}$, where $\mathcal{K} = H^2(\mathbb{D}^k) \otimes H^2(\mathbb{D}^{l-k}) \otimes \underbrace{\mathbb{C}^p \otimes .....\otimes \mathbb{C}^p}_\textrm{m-l times}$.
\end{proposition}
\begin{proof}
By construction, $M=(M_1,....,M_n)$ is a family of $\mathcal{U}_n$-twisted power partial isometries on $\mathcal{H}$ with respect to $\{I_{\mathcal{K}} \otimes U_{ij}\}_{i<j}$. This can be proved by repeated applications of Lemma \ref{mz}, and Lemma \ref{jpd}. For instance, if $1<i<j$, then
$$ M_i^*M_j = (I_{\mathcal{K}} \otimes U_{ji})M_jM_i^*$$ and 

$$M_iM_j = (I_{\mathcal{K}} \otimes U_{ij})M_jM_i.$$ 
\end{proof}

\section{Decomposition of irreducible families of $\mathcal{U}_n$-twisted power partial isometries }
The principal goal of this section is to prove the structure theorem for irreducible families of $\mathcal{U}_n$-twisted power partial isometries. This section is natural continuation of \cite{OPY22}.

\begin{theorem}\label{irred}
Let $\{V_i:i=1,...,n\}$  be irreducible family of $\mathcal{U}_n$-twisted power partial isometries. Consider $V_1$ has Halmos and Wallen decomposition as in Theorem \ref{halwal}.
\begin{enumerate}
\item If $\mathcal{H}=\mathcal{H}_s$, then
$$
V_1=S\otimes 1_{\mathcal{M}_s}, ~~ V_j=D_j(U_{1j})\otimes \tilde{V_j},~~j=2,...,n.
$$
\item If $\mathcal{H}=\mathcal{H}_b$, then
$$
V_1=S^*\otimes 1_{\mathcal{M}_b}, ~~ V_j=D_j(U_{j1})\otimes \tilde{V_j},~~j=2,...,n.
$$
\item If $\mathcal{H}=\mathcal{H}_p$, then
$$
V_1=J_p\otimes 1_{\mathcal{M}_p}, ~~ V_j=d_j(U_{1j})\otimes \tilde{V_j},~~j=2,...,n.
$$
In all cases, $\{\tilde{V_j}, j=2,3,...,n\}$ are irreducible family of $\mathcal{U}_n$-twisted power partial isometries acting on corresponding Hilbert space with $n-1$ generators.
\item If $\mathcal{H}=\mathcal{H}_u$, then
$$
V_1=U_1, ~~ U_1^*V_j=U_{1j}V_jU_1^*,~~ V_jU_1=U_{1j}U_1V_j,~~ j=2,...,n,
$$
where $U_1$ is unitary, operators $V_j$, $j=2,...,n$, are $\mathcal{U}_n$-twisted power partial isometries and the family $\{U_1,V_j,j=2,...,n\}$ is irreducible.
\end{enumerate}
Families corresponding to different cases are non-equivalent. Families corresponding to $\{\tilde{V_j}^{(1)}\}$ and $\{\tilde{V_j}^{(2)}\}$ inside the same case are equivalent if and only if the latter families are equivalent. 

\end{theorem}

\begin{proof}
Assume that $\mathcal{H}=\mathcal{H}_s = \ell^2(\mathbb{Z}_+)\otimes \mathcal{M}_s$ and $V_1 = S\otimes 1_{\mathcal{M}_s}.$ Let $e_n,~~ n\in\mathbb{Z}_+$ be the standard basis of $\ell^2(\mathbb{Z}_+)$. Then we have
$$ \ell^2(\mathbb{Z}_+)\otimes \mathcal{M}_s = \bigoplus_{i=0}^{\infty} e_i \otimes \mathcal{M}_s.$$
Put $\mathcal{H}_i = e_i\otimes \mathcal{M}_s, ~~ i \in \mathbb{Z}_+.$ Then the operators,
$$ P_i = V_1^i (V_1^*)^i - V_1^{i+1}(V_1^*)^{i+1}$$
are orthogonal projections onto $\mathcal{H}_i$. It is easy to observe that 
$$ V_jP_i=P_iV_j,~~~ V_j^*P_i=P_iV_j^*, \qquad i\in \mathbb{Z}_+,~~ j=2,3,...,n.$$
This implies that $\mathcal{H}_i$ is invariant under operators $V_j, V_j^*, j=2,..,n.$ Let $V_j^{(i)}$ denote the restriction of the operator $V_j$ onto $\mathcal{H}_i$. Then we can identify $V_j^{(i)}$ with an operator on $\mathcal{M}_s$  denoted by the same symbol. Then for any $x\in \mathcal{M}_s$, we have
$$ V_jV_1(e_i\otimes x)=e_{i+1}\otimes V_j^{i+1}(x),\qquad V_1V_j(e_i\otimes x)= e_{i+1}\otimes V_j^{(i)}(x).$$
Since $V_jV_1= U_{1j}V_1V_j$, we get $V_j^{(1)}= U_{1j}V_j^{(0)}$, $V_j^{(2)}= U_{1j}V_j^{(1)}= U_{1j}^2V_j^{(0)}$. So, we have $V_j^{(i)}= U_{1j}^iV_j^{(0)}, ~~i\in \mathbb{Z}_+, ~ j=2,...,n$. Put $\tilde{V_j}= V_j^{(0)}$. Then, it is easy to see that
$$V_j = D_j(U_{1j}) \otimes \tilde{V_j},~~~ j=2,...,n.$$

To explore the irreducibility, let $C$ be a operator commuting with\linebreak $V_i, V_i^*,~~i=1,...,n.$ We study a structure of the operator $C$.

In particular, if
$$CV_1=V_1C,~~ CV_1^*=V_1^*C, ~~~ \text{with}~~~ V_1 = S\otimes 1_{\mathcal{M}_s},$$
then, one has
$$ C=1_{\ell^2(\mathbb{Z}_{+})}\otimes \tilde{C}.$$
Then $$CV_j=V_jC,~~ CV_j^*=V_j^*C,~~ j=2,...,n,$$
if and only if
 $$\tilde{C}\tilde{V_j}=\tilde{V_j}\tilde{C},~~ \tilde{C}\tilde{V_j}^*=\tilde{V_j}^*\tilde{C},~~ j=2,...,n.$$
The Schur's lemma implies that $\{V_i:i=1,...,n\}$ is irreducible if and only if $\{\tilde{V_i}:i=2,...,n\}$ is irreducible. Again by similar steps one can show that two families $\{V_i^{(\delta)} :i=1,...,n\}$, $\delta=1,2$, are unitarily equivalent if and only if the corresponding families $\{\tilde{V}_i^{\delta} :i=2,...,n\}$, $\delta=1,2$, are unitarily equivalent.

The remaining cases can also be proved analogously.
\end{proof}

Now we are ready to formulate our classification result.
 
\begin{theorem}
Any irreducible family of $\mathcal{U}_n$-twisted power partial isometries \linebreak $\{V_i:i=1,...,n\}$ is unitarily equivalent to a family of operators acting on 
$$\mathcal{H} = \bigotimes_{i\in \Phi_s}\ell^2(\mathbb{Z}_+)\otimes \bigotimes_{j\in \Phi_b}\ell^2(\mathbb{Z}_+)\otimes \bigotimes_{p,  \Phi_p\neq \emptyset}\mathbb{C}^p \bigotimes\mathcal{H}_u,$$
given by

\begin{align*}
V_j & = \bigotimes_{i<j, i\in\Phi_s}D_j(U_{ij})\otimes S \bigotimes_{i>j, i\in\Phi_s}1_{\ell^2(\mathbb{Z}_+)} \bigotimes_{i\in\Phi_b}1_{\ell^2(\mathbb{Z}_+)}\bigotimes_{i\in \Phi_q,  \Phi_q\neq \emptyset}1_{\mathbb{C}^q} \otimes1_{\mathcal{H}_u}  ~~~  &\text{if}~~ j\in \Phi_s,\\
V_j & = \bigotimes_{i\in\Phi_s}D_j(U_{ij})\bigotimes_{i<j, i\in\Phi_b}D_j(U_{ji})\otimes S^* \bigotimes_{i>j, i\in\Phi_b}1_{\ell^2(\mathbb{Z}_+)}\bigotimes_{i\in \Phi_q,  \Phi_q\neq \emptyset}1_{\mathbb{C}^q} \otimes1_{\mathcal{H}_u} ,  ~~~  &\text{if}~~ j\in \Phi_b,\\
V_j & = \bigotimes_{i\in\Phi_s}D_j(U_{ij})\bigotimes_{i\in\Phi_b}D_j(U_{ji})\bigotimes_{i\in \Phi_q,  \Phi_q\neq \emptyset, i<j}d_j(U_{ij})\otimes J_p\bigotimes_{i\in \Phi_q,  \Phi_q\neq \emptyset, i>j}1_{\mathbb{C}^q} \otimes1_{\mathcal{H}_u}, ~~~ &\text{if}~~ j\in \Phi_p \neq \emptyset,\\
V_j & = \bigotimes_{i\in\Phi_s}D_j(U_{ij})\bigotimes_{i\in\Phi_b}D_j(U_{ji})\bigotimes_{i\in \Phi_q,  \Phi_q\neq \emptyset}d_j(U_{ij}) \otimes U_j,  &\text{if}~~ j\in \Phi_u.
\end{align*}
where
$$ \{1,2,...,n\} = \Phi_s\cup \Phi_b\cup\Psi\cup \Phi_u,$$
the components are disjoint sets and 
$$\Psi = \bigcup_{p=1}^{\infty} \Phi_p$$
with the finite number of non-empty components. The unitary operators \linebreak $\{U_j, ~~j\in \Phi_u\}$ is an irreducible family of unitary operators on $\mathcal{H}_u$ satisfying $U_i^*U_j= U_{ij}U_jU_i^*$, $i\neq j,~~~ i,j \in \Phi_u.$
\end{theorem} 

The following result discusses the uniqueness of the decomposition.

\begin{theorem}
Any irreducible family of $\mathcal{U}_n$-twisted power partial isometries $\{V_i:~ i=1,...,n\}$ is unitarily equivalent to the family described in the above theorem, corresponding to certain decomposition  
$$ \{1,2,...,n\} = \Phi_s\cup \Phi_b\cup \bigcup_{p=1}^{\infty} \Phi_p\cup \Phi_u$$
Families corresponding to different decompositions are non-equivalent. Families corresponding to the same decomposition are equivalent if the related families $\{U_i:~i\in\Phi_u\}$ are equivalent.
\end{theorem}

\begin{remark}
For $i,j\in \{1,2,...,n\}$ with $i\neq j$, suppose $U_{ij}=\lambda_{ij}$ for $\lambda_{ij}\in \mathbb{T}$. Then results in this section imply the main results in the paper \cite{OPY22}.
\end{remark}

\section{Decomposition of $\mathcal{U}_n$-twisted power partial isometries}
In this section, we prove the structure theorem for n-tuple of $\mathcal{U}_n$-twisted power partial isometries without irreducibility assumption.

\begin{theorem}
Let $V=(V_1,...,V_N)$ be a N-tuple of  $\mathcal{U}_N$-twisted power partial isometries on a Hilbert space $\mathcal{H}$, and set $I=\{u,s,b\}\cup \{p\in \mathbb{N}:p\geq 1\}$. For each multiindex $i\in I^N$, let $\Phi_{i,u}:=\{n:1\leq n\leq N, i_n\neq u\}$. Set $\mathcal{K}_{i,n}=\ell^2$ if $i_n=s$ or $i_n=b$, and $\mathcal{K}_{i,n}=\mathbb{C}^{p}$ if $i_n=p$. Also set $D_{\mathcal{K}_{i,m}} = D_j(U_{mn})$ if $i_n=s$, $D_{\mathcal{K}_{i,m}} = D_j(U_{nm})$ if $i_n=b$ and $D_{\mathcal{K}_{i,m}} = d_j(U_{mn})$ if $i_n=p$. Then there are closed subspaces $\{\mathcal{H}_i:i\in I^N\}$ of $\mathcal{H}$. $\{\mathcal{H}_i:i\in I^N\}$ are reducing for $V_n,~\forall~1\leq n\leq N$ and $\mathcal{H}=\bigoplus_{i\in I^N}\mathcal{H}_i$, Hilbert spaces $\{\mathcal{M}_i : i \in I^N\}$, and commuting unitaries $\{T_{i,n}\in U(\mathcal{M}_i): i_n = u\}$ such that the $V_n|_{\mathcal{H}_i}$ for $1\leq n\leq N$ are simultaneously unitary equivalent to
\begin{enumerate}

\item[(a)] $\left(\otimes_{m\in\Phi_{i,u},m< n} D_{\mathcal{K}_{i,m}}\otimes_{m\in\Phi_{i,u},m> n} 1_{\mathcal{K}_{i,m}}\right)\otimes T_{i,n}$ if $i_n = u$;
\item[(b)] $\left(\otimes_{m\in\Phi_{i,u},m< n} D_{\mathcal{K}_{i,m}}\otimes_{m\in\Phi_{i,u},m> n} 1_{\mathcal{K}_{i,m}}\right)\otimes S \otimes 1_{\mathcal{M}_i}$ if $i_n = s$;
\item[(c)] $\left(\otimes_{m\in\Phi_{i,u},m< n} D_{\mathcal{K}_{i,m}}\otimes_{m\in\Phi_{i,u},m> n} 1_{\mathcal{K}_{i,m}}\right)\otimes S^* \otimes 1_{\mathcal{M}_i}$ if $i_n = b$;
\item[(d)] $\left(\otimes_{m\in\Phi_{i,u},m< n} D_{\mathcal{K}_{i,m}}\otimes_{m\in\Phi_{i,u},m> n} 1_{\mathcal{K}_{i,m}}\right)\otimes J_p \otimes 1_{\mathcal{M}_i}$ if $i_n = p$.
\end{enumerate}
\end{theorem}
\begin{proof}
We will prove this theorem by induction on $N$. Assume that the theorem holds for $N$-tuples of $\mathcal{U}_N$-twisted power partial isometries and the subspaces $\mathcal{H}_i$ are reducing for every operator $W$ that is $\mathcal{U}_{N+1}$-twisted with all the \linebreak $V_n, ~1\leq n\leq N$. For $N=1$, the sets $i$ are singletons. The subspaces $\mathcal{H}_i$ are the subspaces $\mathcal{H}_u$, $\mathcal{H}_s$, $\mathcal{H}_b$ and $\mathcal{H}_p,~p\geq 1$  as in Theorem $\ref{halwal}$. Suppose  $W$ is $\mathcal{U}_2$-twisted with $V$ then by Lemma  \ref{reducing} the subspaces $\mathcal{H}_u$, $\mathcal{H}_s$, $\mathcal{H}_b$ and $\mathcal{H}_p,~p\geq 1$ are reducing for $W$.

For $N=2$, let $(V_1,V_2)$ be a pair of $\mathcal{U}_2$-twisted power partial isometries on a Hilbert space $\mathcal{H}$, and set $I=\{u,s,b\}\cup \{p\in \mathbb{N}:p\geq 1\}$. Since $V_1$ is a power partial isometry by Theorem \ref{halwal}, we have 

\begin{equation*}
\begin{split}
\mathcal{H} &= \mathcal{H}_u \oplus \mathcal{H}_s \oplus \mathcal{H}_b \oplus \left(\bigoplus_{p=1}^{\infty} \mathcal{H}_p\right)\\
&=\mathcal{M}_u \oplus (\ell^2\otimes \mathcal{M}_s) \oplus (\ell^2\otimes \mathcal{M}_b) \oplus \left(\bigoplus_{p=1}^{\infty} \mathbb{C}^p \otimes \mathcal{M}_p\right),\\
\end{split}
\end{equation*}
where $\mathcal{M}_u=\mathcal{H}_u$, and 
$$V_1 = T \oplus (S\otimes 1_{\mathcal{M}_s}) \oplus (S^*\otimes 1_{\mathcal{M}_b}) \oplus \left(\bigoplus_{p=1}^{\infty} J_p\otimes 1_{\mathcal{M}_p}\right)$$
 where $T$ is a unitary operator on $\mathcal{M}_u$.
Now, since $\mathcal{M}_i, i\in I$ reduces $V_2$, again applying Theorem \ref{halwal} we yield
$$\mathcal{M}_i = \mathcal{N}_{(i,u)}\oplus (\ell^2\otimes\mathcal{N}_{(i,s)})\oplus (\ell^2\otimes\mathcal{N}_{(i,b)})\oplus\left(\bigoplus_{p=1}^{\infty} (\mathbb{C}^p\otimes\mathcal{N}_{(i,p)})\right).$$
Therefore 
$$\mathcal{H}=\bigoplus_{i\in I^2}\mathcal{H}_i =\bigoplus_{i\in I^2}\left(\otimes_{m\in\Phi_{i,u}}\mathcal{K}_{i,m}\right)\otimes \mathcal{N}_i.$$ 
Suppose $V_1|_{\mathcal{H}_{i}}$ is a shift, then we have 
$$\mathcal{H}_i = \ell^2 \otimes  \mathcal{M}_i,$$
and 
$$ V_1|_{\mathcal{H}_{i}} = S \otimes 1_{\mathcal{M}_i}.$$ 
From Theorem \ref{irred}  
$$V_2|_{\mathcal{H}_{i}} = D_j(U) \otimes \tilde{V_2}.$$
Similarly, if $V_1|_{\mathcal{H}_{i}}$ is a backward shift, then we have 
$$\mathcal{H}_i = \ell^2 \otimes  \mathcal{M}_i,$$
and 
$$ V_1|_{\mathcal{H}_{i}} = S^* \otimes 1_{\mathcal{M}_i}.$$ 
From Theorem \ref{irred}
$$V_2|_{\mathcal{H}_{i}} = D_j(U^*) \otimes \tilde{V_2}.$$
Now, if $V_1|_{\mathcal{H}_{i}}$ is a truncated shift, then we have 
$$\mathcal{H}_i = \mathbb{C}^p \otimes  \mathcal{M}_i,$$
and 
$$ V_1|_{\mathcal{H}_{i}} = J_p \otimes 1_{\mathcal{M}_i}.$$ 
Again from Theorem \ref{irred} 
$$V_2|_{\mathcal{H}_{i}} = d_j(U) \otimes \tilde{V_2}.$$
Then 
$$V_1= \begin{cases} 
      \bigoplus_{i \in I^2}\left( 1_{\mathcal{K}_{i,i_2\neq u}}\otimes T_{i,1}\right) & \text{if}~ i_1 = u, \\
      \bigoplus_{i \in I^2}\left(1_{\mathcal{K}_{i,i_2\neq u}}\otimes S\otimes 1_{\mathcal{N}_i}\right) &  \text{if}~ i_1 = s,\\
      \bigoplus_{i \in I^2}\left(1_{\mathcal{K}_{i,i_2\neq u}}\otimes S^*\otimes 1_{\mathcal{N}_i}\right) &  \text{if}~ i_1 = b,\\
      \bigoplus_{i \in I^2}\left(1_{\mathcal{K}_{i,i_2\neq u}}\otimes J_p\otimes 1_{\mathcal{N}_i}\right) &  \text{if}~ i_1 = p
   \end{cases}
$$
and 
$$V_2= \begin{cases} 
      \bigoplus_{i \in I^2}\left( D_{\mathcal{K}_{i,i_1\neq u}}\otimes T_{i,2}\right) & \text{if}~ i_2 = u \\
      \bigoplus_{i \in I^2}\left(D_{\mathcal{K}_{i,i_1\neq u}}\otimes S\otimes 1_{\mathcal{N}_i}\right) &  \text{if}~ i_2 = s,\\
      \bigoplus_{i \in I^2}\left(D_{\mathcal{K}_{i,i_1\neq u}}\otimes S^*\otimes 1_{\mathcal{N}_i}\right) &  \text{if}~ i_2 = b,\\
      \bigoplus_{i \in I^2}\left(D_{\mathcal{K}_{i,i_1\neq u}}\otimes J_p\otimes 1_{\mathcal{N}_i}\right) &  \text{if}~ i_2 = p.
   \end{cases}
$$
Suppose $W$ is $\mathcal{U}_3$-twisted with $(V_1,V_2)$. Then by Lemma \ref{reducing} and induction hypothesis, each $\mathcal{H}_i$ are reducing for $W$.

Suppose the theorem is true for $N$. We will prove the theorem for $N+1$. Let $V=(V_1,...,V_{N+1})$ be a $\mathcal{U}_{N+1}$-twisted power partial isometries. We set \linebreak $S_s:=S,~S_b:=S^*$ and $S_p:=J_p$ for $p\geq 1$ for the simplification of notations and reducing the number of cases. We apply the induction hypothesis to $V=(V_1,...,V_{N})$. Since we are assuming the simultaneous unitary equivalence. We can conjugate the operators by a unitary. Suppose there exist Hilbert spaces $\mathcal{N}_i$ with

$$\mathcal{H}=\bigoplus_{i\in I^N}\mathcal{H}_i =\bigoplus_{i\in I^N}\left(\bigotimes_{m\in\Phi_{i,u}}\mathcal{K}_{i,m}\right)\otimes \mathcal{N}_i,$$ 
and that for every $1 \leq n \leq N$,
$$V_n= \begin{cases} 
      \bigoplus_{i \in I^N}\left(\otimes_{m\in\Phi_{i,u},m< n} D_{\mathcal{K}_{i,m}}\otimes_{m\in\Phi_{i,u},m> n} 1_{\mathcal{K}_{i,m}}\right)\otimes T_{i,n} & \text{if}~ i_n = u \\
      \bigoplus_{i \in I^N}\left(\otimes_{m\in\Phi_{i,u},m< n} D_{\mathcal{K}_{i,m}}\otimes_{m\in\Phi_{i,u},m> n} 1_{\mathcal{K}_{i,m}}\right)\otimes S \otimes 1_{\mathcal{M}_i}& \text{if}~ i_n = s,\\
      \bigoplus_{i \in I^N}\left(\otimes_{m\in\Phi_{i,u},m< n} D_{\mathcal{K}_{i,m}}\otimes_{m\in\Phi_{i,u},m> n} 1_{\mathcal{K}_{i,m}}\right)\otimes S^* \otimes 1_{\mathcal{M}_i}& \text{if}~ i_n = b,\\
      \bigoplus_{i \in I^N}\left(\otimes_{m\in\Phi_{i,u},m< n} D_{\mathcal{K}_{i,m}}\otimes_{m\in\Phi_{i,u},m> n} 1_{\mathcal{K}_{i,m}}\right)\otimes J_p \otimes 1_{\mathcal{M}_i}& \text{if}~ i_n = p.
   \end{cases}
$$

The power partial isometry $V_{N+1}$ is $\mathcal{U}_{N+1}$-twisted with $V_n$ for $1\leq n\leq N$. Using Lemma \ref{reducing} and induction hypothesis, all the summands $\mathcal{H}_i$ are reducing for $V_{N+1}$. Therefore $V_{N+1}|_{\mathcal{H}_i}$ is $\mathcal{U}_{N+1}$-twisted with all the operators of the form
 $\left(\otimes_{m\in\Phi_{i,i_n},m< n} D_{\mathcal{K}_{i,m}}\otimes_{m\in\Phi_{i,i_n},m> n} 1_{\mathcal{K}_{i,m}}\right)\otimes S_{i_n}\otimes 1_{\mathcal{N}_i}$
coming from the summands of the $V_n|_{\mathcal{H}_i}$. Let $V$ be in the $C^*$-subalgebra of $B\left(\bigotimes_{j\in\Phi_i}\mathcal{K}_{i,j}\right)$ generated by the operators of the form 
$\left(\otimes_{m\in\Phi_{i,i_n},m< n} D_{\mathcal{K}_{i,m}}\otimes_{m\in\Phi_{i,i_n},m> n} 1_{\mathcal{K}_{i,m}}\right)\otimes S_{i_n}$. Thus $V_{N+1}|_{\mathcal{H}_i}$ is $\mathcal{U}_{N+1}$-twisted with all the operators of the form $V\otimes 1$.

If $i_n=p\geq 1$, the $C^*$-algebra generated by $S_{i_n}$ is $C^*(S_{i_n})=C^*(J_p)$, thus $C^*(S_{i_n})=M_p(\mathbb{C})$. If $i_n=s$ or $i_n=b$, $C^*(S_{i_n})=C^*(S)$ contains the compact operators on $\ell^2$. Thus the $C^*$-algebra acts irreducibly on $\mathcal{K}_{i,n}$ for all $n\in \Phi_{i,u}$. Therefore the spatial tensor product of $C^*$-algebras $\bigotimes_{n\in \Phi_{i,u}}C^*(S_{i_n})$ acts irreducibly on $\bigotimes_{n\in \Phi_{i,u}} \mathcal{K}_{i,n} $. Hence the operator $V_{N+1}|_{H_i}$ has the form $\bigotimes_{n\in\Phi_{i,u}}D_{\mathcal{K}_{i,n}}\otimes R_i$ for some $R_i\in B(\mathcal{N}_i)$.  Since $\mathcal{H}_i$ is a reducing subspace for $V_{N+1}$, we have $V_{N+1}|_{\mathcal{H}_i}$  is a power partial isometry. Thus $R_i$ is a power partial isometry. Since the $V_n$ is $\mathcal{U}_{N+1}$-twisted with $V_{N+1}$, thus the unitaries $T_{i,n}$ is $\mathcal{U}_{N+1}$-twisted with $R_i$.

Since $R_i$ is a power partial isometry on $\mathcal{N}_i$. By Theorem \ref{halwal}, we get a direct sum decomposition of $\mathcal{N}_i$. Again we can conjugate by a unitary isomorphism. Let $\mathcal{N}_i$ has a decomposition as follows:

$$\mathcal{N}_i = \mathcal{M}_{(i,u)}\oplus (\ell^2\otimes\mathcal{M}_{(i,s)})\oplus (\ell^2\otimes\mathcal{M}_{(i,b)})\oplus\left(\bigoplus_{p=1}^{\infty} (\mathbb{C}^p\otimes\mathcal{M}_{(i,p)})\right)~ \text{and}$$
\begin{equation}\label{eqsplit1}
R_i=U_{i,u} \oplus (S\otimes D_{\mathcal{M}_{(i,s)}})\oplus (S^*\otimes D_{\mathcal{M}_{(i,b)}})\oplus \left(\bigoplus_{p=1}^{\infty} (J_p\otimes D_{\mathcal{M}_{(i,p)}})\right)
\end{equation}

with $U_{i,u}=R_i|_{\mathcal{M}_{(i,u)}}$ unitary. Now for $i'=(i,i_{N+1})\in I^N \times I = I^{N+1}$, we take

\begin{equation*}
\mathcal{H}_{i'}= \begin{cases} 
      \mathcal{H}_i=\left(\bigotimes_{n\in\Phi_{i,u}}\mathcal{K}_{i,n}\right)\otimes \mathcal{M}_{(i,u)}  & \text{if}~ i_{N+1} = u \\
    \left(\bigotimes_{n\in \Phi_{i,i_{N+1}}} {\mathcal{K}_{i,n}}\right)\otimes \mathcal{K}_{i',N+1}\otimes \mathcal{M}_{(i,i_{N+1})} &  \text{if}~ i_{N+1} \neq u.
   \end{cases}
   \end{equation*}
 Since then $\Phi_{i',u}=\Phi_{i,u}$, we have  
 \begin{equation}\label{eqsplit2}
 \mathcal{H}_{i'} =\begin{cases} 
      \left(\bigotimes_{n\in\Phi_{i',u}}\mathcal{K}_{i',n}\right)\otimes \mathcal{M}_{(i,u)}  & \text{if}~ i_{N+1} = u  \\
    \left(\bigotimes_{n\in \Phi_{i',i_{N+1}}} {\mathcal{K}_{i',n}}\right)\otimes \mathcal{M}_{(i,i_{N+1})} &  \text{if}~ i_{N+1} \neq u.
   \end{cases}
   \end{equation}
From equations \ref{eqsplit1} and \ref{eqsplit2}, we get $\mathcal{H}_i=\bigoplus_{j\in I }\mathcal{H}_{(i,j)}$ for each $i\in I^N$, thus $\mathcal{H}=\bigoplus_{i'\in I^{N+1} }\mathcal{H}_{i'}.$ Since $T_{i,n}$ is $\mathcal{U}_{N+1}$-twisted with $R_i$ for all $i\in I^N$. From the case $N=1$, observe that $T_{1,n}$ is $\mathcal{U}_{2}$-twisted $U_{1,u}$. Thus all the direct summands of subspaces in \ref{eqsplit1} are reducing for $T_{1,n}$. If $i_{N+1}=u$, then we take $T_{i',m}$ to be $T_{i,m}$. If $i_{M+1}\neq u$, then we take $T_{i',m}$ to be the operator on $\mathcal{M}_{i'}=\mathcal{M}_{(i,i_{N+1})}$ such that $T_{i,n}|_{\mathcal{K}_{i,n}\otimes\mathcal{M}_{(i,i_{N+1})}}=D_{\mathcal{K}_{i,n}}\otimes T_{i',m}$.

It is necessary to check that, if any operator $W$ is $\mathcal{U}_{N+2}$-twisted with \linebreak $\{V_n:1\leq n\leq N+1\}$ then the  subspaces $\{\mathcal{H}_i:i\in I^{N+1}\}$ are reducing for the operator $W$. By Lemma \ref{reducing} and inductive hypothesis, observe that each $\mathcal{H}_i$ are reducing for $W$. Thus $W|_{\mathcal{H}_i}$ has the form $\bigotimes_{n\in\Phi_{i,u}}D_{\mathcal{K}_{i,n}}\otimes R$ for some $R\in B(\mathcal{N}_i)$. Since $W|_{\mathcal{H}_i}=\bigotimes_{n\in\Phi_{i,u}}D_{\mathcal{K}_{i,n}}\otimes R$ is $\mathcal{U}_{N+2}$ twisted with \linebreak $V_{N+1}|_{H_i}=\bigotimes_{n\in\Phi_{i,u}}D_{\mathcal{K}_{i,n}}\otimes R_i$. From Lemma \ref{reducing} and the case $N=1$, it follows that the subspaces in the decomposition \ref{eqsplit1} are all reducing for $R$. This proves our induction hypothesis for $N+1$. Hence proof completes.
\end{proof}

{\bf Acknowledgments.} We would like to thank Jaydeb Sarkar for suggesting the problem and for valuable discussions. The first author is supported by the Junior Research Fellowship (09/0239(13298)/2022-EM) of CSIR (Council of Scientific and Industrial Research, India). The second author is supported by the Teachers Association for Research Excellence (TAR/2022/000063) of SERB (Science and Engineering Research Board, India).

\nocite{*}
\bibliographystyle{amsplain}
\bibliography{database}  

\end{document}